\documentclass[10pt,reqno]{amsart}
\usepackage[cmtip,arrow]{xy}
\usepackage{pb-diagram}
\usepackage[dvips]{graphics}
\usepackage[dvips]{epsfig}
\usepackage[latin1]{inputenc}
\usepackage{amsmath}
\usepackage{amssymb}
\usepackage{enumitem}

\newtheorem{theorem}{Theorem}[section]

\newtheorem{lemma}[theorem]{Lemma}
\newtheorem{prop}[theorem]{Proposition}

\newtheorem{remark}[theorem]{Remark}
\newtheorem{corollary of proof}[theorem]{Corollary of Proof}

\newcommand{\Gl}{\operatorname{Gl}}
\newcommand{\gl}{\mathfrak{gl}}

\newcommand{\R}{\mathbb{R}}
\newcommand{\C}{\mathbb{C}}

\newcommand{\N}{\mathbb{N}}
\newcommand{\Z}{\mathbb{Z}}

\renewcommand{\qed}{$\hfill \square$ \medskip \\}

\newcommand{\piinv}{\pi^{-1}}

\renewcommand{\proof}{{\bf Proof: }}

\newcommand{\ad}{\operatorname{ad}}
\newcommand{\tr}{\operatorname{tr}}

\newcommand{\id}{\operatorname{id}}
\newcommand{\g}{\mathfrak{g}}

\newcommand{\nach}{\to}

\newcommand{\dt}{\frac{d}{dt}}
\newcommand{\dtnull}{\left. \dt \right|_{t=0}}
\newcommand{\dttnull}{\left. \dt \right|_{t=t_0}}

\newcommand{\tensor}{\otimes}
\newcommand{\omA}{\omega_A}

\renewcommand{\phi}{\varphi}

\newcommand{\Hom}{\operatorname{Hom}}

\newcommand{\tildeu}{\tilde{u}}
\newcommand{\tildew}{\tilde{w}}
\newcommand{\tildev}{\tilde{v}}

\newcommand{\dnablaA}{d^{\nabla^{A}}}
\newcommand{\nablaA}{\nabla^{A}}

\renewcommand{\t}{h_{\langle \, . , . \, \rangle}}

\renewcommand{\k}{\mathfrak{k}}

\newcommand{\Ja}{J_\alpha}
\newcommand{\Njalpha}{N_{J_\alpha}}
\newcommand{\ga}{g_\alpha} 
\newcommand{\val}{^\vee \alpha}
\renewcommand{\v}{^\vee}
\newcommand{\w}{^\wedge}
\newcommand{\cdAh}{\v \circ d_A \circ \w}

\newcommand{\Psom}{P_{SO(M)}}
\newcommand{\phitJaX}{\phi^t_{\Ja X^\#}}
\newcommand{\JaX}{\Ja X^\#}
\newcommand{\cdt}{\dot{c}(t)}

\newcommand{\zeth}{ [\zeta]_{dR}(H_1(Y;\Z)) }
\newcommand{\clzeth}{ \overline{\zeth}}
\newcommand{\Cquotzetah}{\C /{( \zeth})}
\newcommand{\Cquotclzetah}{\C / {(\overline{\zeth})}}

\renewcommand{\Re}{\mathfrak{Re}}
\renewcommand{\Im}{\mathfrak{Im}}

\newcommand{\so}{\mathfrak{so}}
\newcommand{\su}{\mathfrak{su}}

\begin{document}
\thispagestyle{empty}
\title[Integrable a.c. structures in principal bundles]{Integrable almost complex structures in principal bundles and holomorphic curves} 
\author[Raphael Zentner]{Raphael Zentner} 
%\date{Started 9th March 2010}
%\address {Mathematisches Institut \\ Westf\"alische Wilhelms-Universit\"at \\ Einsteinstr. 62 \\ 48149 M\"unster\\
%Germany}

%\email{raphael.zentner@math.uni-muenster.de}
\date{June 2003, new layout August 2010, revision November 2011, 2nd revision April 2012}

\maketitle

\begin{abstract}
We consider almost complex structures that arise naturally in a particular class of principal fibre bundles, where the choice of a connection can be used to determine equivariant isomorphisms between the vertical and horizontal tangent bundles of the total space. For instance, such data always exist on the frame bundle of a 3-manifold, but also in many other situations. We study the integrability condition to a complex structure, obtaining a system of gauge invariant coupled first order partial differential equations. This yields to a few correspondences between complex-geometric properties on the total space and metric properties on the base.
\end{abstract}

\section*{Introduction}
\thispagestyle{empty}
In this article, we consider almost complex structures that arise naturally in a particular class of principal fibre bundles, where the choice of a connection can be used to determine equivariant isomorphisms between the vertical and horizontal tangent bundles of the total space. More precisely, we assume given a $K$-principal fibre bundle $\pi: P \to M$ with a connection $A$ on $P$, where $K$ is a compact Lie group. Denoting by $\ad(P)$ the adjoint bundle, we suppose given an $\ad(P)$-valued 1-form $\alpha$ which defines a bundle isomorphism $TM \nach \ad(P)$. To such a triple $(\pi, A, \alpha)$ we will associate a $K$-invariant almost complex structure $\Ja$ on the total space $P$. A necessary condition for the existence of such a triple is that $P$ is (almost complex) parallelisable. We will study the integrability of this almost complex structure. This will yield a system of gauge invariant coupled non-linear partial differential equation for $(\alpha,A)$ in which the curvature of $A$ and the Lie algebra structure of $K$ appear. 

The necessary topological setup, yielding an isomorphism $\alpha$ as needed in our construction, occurs for instance for $P$ the oriented orthonormal frame bundle of a $3$-manifold. Also, if $P$ is the frame bundle of a sum of $n$ complex line bundles, where $n$ is the dimension of a parallelisable $n$-dimensional manifold, with the natural reduction of its structure group to $T^n$, then we are in such a situation. Another class of exmples comes from oriented Riemannian $6$-manifolds that have a reduction of its structure group from $SO(6)$ to $SO(3) \times SO(3)$. 

Via $\alpha$, an $\ad$-invariant inner product on the Lie algebra $\k$ of $K$ induces the structure of an Euclidean vector bundle on $\ad(P)$, and therefore we get an induced metric $g_\alpha$ on the base $M$. We show that the pull-back via $\alpha$ of the linear connection on $\ad(P)$, induced from $A$, is the Levi-Civita-Connection of $(M,g_\alpha)$ if and only if the first of our integrability equations is satisfied. The second integrability equation implies that the sectional curvature of $g_\alpha$ is non-positive in general. In the case $K= SU(2)$ or $K=SO(3)$, the integrability equations are both satisfied if and only if the pull-back connection is the Levi-Civita connection and the metric $g_\alpha$ is hyperbolic. Therefore we get a  relation between our integrability problem and hyperbolic geometry. In particular, we recover the fact that the bundle of orthonormal oriented frames of every hyperbolic oriented 3-manifold is naturally a complex manifold. We also show that there is an equivalence between complete integrability of $\Ja$-holomorphic vector fields and the geodesic completeness of the base manifold, assuming only that the first of our integrability equations holds.

Assuming geodesic completeness of the base and that the integrability conditions hold, we obtain a holomorphic action of the complexified Lie group $G= K^\C$ on the complex manifold $P$. This action is locally free and transitive if $M$ is connected. By a result of Wang \cite{Wan}, any closed complex manifold with complex parallelisable tangent bundle is in fact given by the quotient of a complex Lie group by a discrete subgroup. However, the proof relies essentially on the compactness of the complex manifold, implying that certain complex-valued holomorphic functions on the manifold are necessarily constant by Liouville's theorem. Our result here may be seen as a generalisation of this result in the particular setting of an integrable almost complex structure $\Ja$ as studied here. In fact, requiring the base manifold $(M,g_\alpha)$ to be geodesically complete is a much weaker condition than requiring $P$ to be closed. 

Closed complex manifolds that are complex parallelisable have been extensively studied by Winkelmann \cite{Wi1, Wi2, Wi3, Wi4}. Our approach may be seen as an alternative viewpoint. In fact, they are quite closely related -- one may ask how general our construction is in the situation that we have a complex parallelisable complex manifold $Q$ that is given as the quotient of a complex Lie group $K^\C$ by a discrete subgroup, and on which the Lie group $K$ acts freely. We will show that in this situation the natural almost complex structure induced by the complex structure of the manifold $Q$ may be obtained as in our setting, namely, as a $K$-invariant almost complex structure $J_{\bar{\alpha}}$ on the $K$-principal bundle $Q \to Q/K$ associated to a bundle isomorphism $\bar{\alpha}: T(Q/K) \to \ad(Q)$, that itself is induced from a `sample' situation $\alpha: T(K^\C/K) \to \ad_K(K^\C)$.

This description as a homogeneous space also gives a $\g$-valued holomorphic form $\omega$ whose real part is just the connection form of $A$. 
We will finally discuss some restrictions on holomorphic maps from Riemann surfaces $\Sigma$ (or more general complex manifolds $Y$) to $P$. First we show that if $f: \Sigma \to P$ is a holomorphic map such that $f^* \omega$ is of scalar form -- i.e. $f^* \omega = Z \cdot \zeta$ for a constant element $Z \in \g$ and a holomorphic form $\zeta$--, then $f$ factorises through an elliptic curve. Finally we consider holomorphic maps $f: \Sigma \to P$ such that the induced map $\bar{f}: \Sigma \nach M$ is conformal with respect to the canonical almost complex structure $J_0$ of $\Sigma$ (i.e. $ \bar{f}^*g_\alpha \ (u,J_0 u) = 0 $ for all $u \in TY$). This differential geometric hypothesis then has an interesting algebraic geometric consequence: If $f$ is non-constant, then the holomorphic map $\Sigma \nach P(\g)$ associated with $f^*\omega$ factorises through a smooth quadric. 
\\

As a perspective, our setup allows to consider natural generalisations of complex-geometric questions even if the integrability conditions are not satisfied. For instance, they may be used to quantify the `deviation' from integrability. Or one might also consider (pseudo-) holomorphic curves for non-integrable almost complex structures in our setting. 
\\

In the first section will set up the notation and introduce our construction explicitely. We establish the integrability condition of $\Ja$. In the second section, we will set up the relation to the geometry on the base manifold, and we will give a complete solution on the frame bundle of a hyperbolic 3-manifold. There we will also describe other situations with the required topological setup. In the third section we will relate (complex) integrability of almost complex vector fields to geodesic completeness of the base. Assuming the integrability conditions hold, this yields the description of the manifold as the quotient of a complex Lie group by a discrete subgroup. In the fifth section, we shall use these results to establish the mentioned complex geometric properties.

\section*{Acknowledgements}
I am delighted to thank Andrei Teleman for multiple essential ideas (this work arose in the beginning of the author's PhD thesis with him as advisor). Furthermore, I would like to express my gratitude to the very useful comments of an anonymous referee.

\section{A class of almost complex structures on principal fibre bundles}
We shall consider situations, in which a principal bundle $P \to M$ on a manifold $M$ with compact structure group $K$, equipped with a connection $A$, admits admits an almost complex structure of a special type. Namely, the almost complex structure will be compatible with the action of the structure group, and it will interchange vertical and horizontal tangent vectors. In this section, we will make this setting more concrete, and we shall examine the integrability condition of this almost complex structure to a complex structure.

\subsection{Notations and main construction}
If one is given an isomorphism of vector spaces 
$ \alpha: V \stackrel{\cong}{\nach} W $, then 
$$
J_\alpha = \begin{pmatrix} 0 & -\alpha^{-1} \cr \alpha & 0 \end{pmatrix}
$$
defines a complex structure on the vector space $V\oplus W$.
Let $K$ be a compact Lie group (or more generally a Lie group admitting an inner product on its Lie algebra that is invariant under the adjoint action), and suppose given a $K$-principal fibre bundle $\pi: P \to M$, where $P, M$ are supposed to be  differentiable
($C^\infty$) manifolds. We denote the smooth free right-action by $\rho: P\times K \nach P$ . We suppose also given a connection $A$ on $P$, and its connection 1-form on $P$ will be denoted by $\omega_A$. This gives a splitting of its tangent bundle into the horizontal and the vertical subbundle
$$
TP = A \oplus V \ .
$$ 
(By definition $V$ is the kernel bundle of $\pi_*:TP \nach TM$).
Thus if $\dim(K) = \dim(M)$ and if one has given $\forall p \in P$ an isomorphism
$$
\alpha_p : A_p \stackrel{\cong}{\nach} V_p
$$
in a differentiable manner, then one is able to construct an almost complex structure on $P$. In this setting it is natural to suppose that $\alpha$ commutes with the given right-action of $K$ on $P$. 
\\

Let $\k$ denotes the Lie algebra of $K$. To each element $X \in \k$ corresponds a vertical vector field, noted $X^\#$ and these vector fields span the vertical bundle, hence it is trivial. Furthermore $\#$ is a Lie algebra homomorphism. In fact one differentiates the right-action $P \times K \stackrel{\rho}{\nach} P$ in the second argument at the identity to obtain a mapping $P \times \k \nach TP$, that can also be understood as a map $\k \nach \Gamma(P,TP)$. These vector fields are equivariant with respect to right-translations in the following sense
\begin{eqnarray}\label{equiv_dieze}
\begin{diagram}
\node{\k} \arrow[1]{e,t}{\#_{pg}}  \node[1]{V_{pg}}  \\
\node{\k} \arrow[1]{n,l}{\ad_{g^{-1}}} \arrow{e,b}{\#_p} \node[1]{V_p} 
\arrow[1]{n,r}{(R_g)_{*,p}}
\end{diagram} 
\end{eqnarray}
where $R_g(p) = \rho(p,g) \equiv pg$ and $(R_g)_{*,p}$ is its derivation at $p$. 

We denote by $A^k_{\ad}(P,\k)$ the space of `tensorial k-forms of type $\ad$' (see \cite{KN1}). By definition, such a form is zero on vertical vectors in $P$ and that it is of type $\ad$ means that the following diagram commutes for all $p \in P$ and $g \in K$,

\begin{eqnarray}\label{equiv_alpha}
\begin{diagram}
\node{A_{pg}} \arrow[1]{e,t}{\alpha_{pg}}  \node[1]{\k}  \\
\node{A_p} \arrow[1]{n,l}{(R_g)_{*,p}} \arrow{e,b}{\alpha_p} \node[1]{\k \, .} 
\arrow[1]{n,r}{\ad_{g^{-1}}} 
\end{diagram}
\end{eqnarray}
\\

In the sequel we assume having a form $\alpha \in A^1_{\ad}(P,\k)$ such that 
$$
\alpha_p|_{A_p} : A_p \nach \k
$$
is an isomorphism for each $p \in P$. Of course, the existence of such an $\alpha$ is a strong topological hypothesis. Indeed this means that the vector bundles $TM$ and the associated bundle $\ad(P) = P \times_{\ad} \k$, the bundle associated to the adjoint representation of $K$ on its Lie algebra, must be isomorphic. Then $P$ is also parallelisable.
The compostion $\# \circ \alpha \equiv \alpha^\#$  defines a $K$-invariant bundle isomorphism $A \to V$ with the required properties (we denote $\alpha$ instead of $\alpha|_{A}$, but we bear in mind that even though $\alpha$ is defined without specifying a connection, we need to have specified a connection so that $\alpha^{-1}$ makes sense). We may thus define
\begin{eqnarray*}
        (J_\alpha)_p : \left\{ \begin{array}{ccc} T_pP & \longrightarrow & T_pP  \\
               \left( \begin{array}{c} w_p \\ X^\#_p \end{array} \right) &
\longmapsto & \left( \begin{array}{cc} 0 & - \alpha^{-1}_p \\ (\alpha^\#)_p & 0 \end{array} \right)   \left( \begin{array}{c} w_p \\ X \end{array}
\right) \ .
        \end{array} \right.
\end{eqnarray*}
$J_\alpha$ is an almost-complex structure on $TP$. 
\subsection{Integrability condition}
The most natural question to ask from the viewpoint of complex geometry is whether $J_\alpha$ is integrable, this means, whether $J_\alpha$ comes from the structure of a complex manifold on $P$.
\\

We review some standard notation. The curvature $F_A$ of a connection $A$ with connection $1$-form $\omega_A$ is given by the formula
\[
 F_A = d \omega_A + \frac{1}{2} \, [\omega_A \wedge \omega_A].
\]
In this formula, we have used the `hybrid notation' $ [\, . \, \wedge \, . \, ]$ denoting the tensor product of wedge-product on the $1-$form factors and the Lie bracket on the Lie algebra factor. We also recall that the differential $d_A$ associated with the curvature $d_A$ is given by $\pi_A^* d$, where $\pi_A$ denotes the projection $TP \to A$ along the vertical bundle $V$.

We can now state our integrability condition of $\Ja$ in terms of the Lie algebra valued 1-form $\alpha$ and the connection $A$. 
\begin{theorem}\label{main_prop}
The almost-complex structure $J_\alpha$ is integrable to a complex structure if and only if the following two equations hold:
\begin{equation} \label{integrability equations}
	\begin{split}
        d_A \alpha & =  0 \\
        F_A & =  \frac{1}{2} \, [\alpha \wedge \alpha] \ .
	\end{split}
\end{equation}
\end{theorem}
\begin{remark}
It is worth pointing out that if the integrability condition is satisfied the curvature of the connection $A$ is entirely determined by $\alpha$ and the Lie algebra $\k$.
\end{remark}

\proof
We recall that $J_\alpha$ is integrable if and only if the Nijenhuis tensor $N_{J_\alpha}$, defined by
$$
   \Njalpha (\zeta, \xi)= [\Ja \zeta, \Ja \xi ] - [\zeta, \xi] - \Ja [\zeta, \Ja \xi] 
        - \Ja [ \Ja \zeta, \xi ]
$$
for  $\zeta, \xi \in \Gamma(P,TP)$, vanishes identically. This is the famous theorem of Newlander and Nirenberg \cite{NN}. We prove the Theorem by showing that the Nijenhuis tensor vanishes if and only if the equations (\ref{integrability equations}) hold. 

As $\Njalpha$ is a tensor, its value $ \Njalpha (\zeta, \xi)_p$ at a given point $p \in P$ does in fact only depend on the elements $\zeta_p, \xi_p \in T_pP$, and not on the vector fields $\zeta, \xi$ extending these vectors in $T_pP$. More precisely, if $\overline{\zeta}$ and $\overline{\xi}$ are other vector fields on $P$ with $\overline{\zeta}_p = \zeta_p$ and $\overline{\xi}_p = \xi_p$, then $\Njalpha (\zeta, \xi)_p = \Njalpha (\overline{\zeta},\overline{\xi})_p$. We will use this fact by choosing vector fields that are particularly suitable to our setting, namely the fundamental vector fields of the type $X^\#$ with $X \in \k$. Furthermore, the Nijenhuis tensor vanishes if and only if 
\begin{equation*}
	\Njalpha (X^\#,Y^\#) = 0 
\end{equation*}
everywhere, for all elements $X, Y$ of the Lie algebra of $K$. In fact, at each point $p \in P$ the vectors of the form $X^\#_p$ and $\Ja X^\#_p$ span the tangent space $T_pP$. Furthermore, the (antisymmetric) Nijenhuis tensor $N_J$ of an almost complex structure $J$ satisfies the identity
\begin{equation*}
	 N_J (\zeta, J \xi) = -J \, N_J (\zeta,\xi) \ 
\end{equation*} 
for all vector fields $\zeta, \xi$.\\

%
%
% We denote $\pi_A$ respectively $\pi_V$ the projection from $TP$ to the horizontal subbundle $A$ respectively to the vertical subbundle $V$. We notify that the vertical projection is just given in terms of the connection form by
%$$
%\pi_V = \# \circ \omega_A 
%$$
%or more precisely $\pi_V (\zeta_p) = (\omega_A (\zeta_p))^\#_p$ for $\zeta_p \in T_pP$.
%%
%%
%%
%%
%%
%%

Applying the definition of the almost complex structure $\Ja$ yields
\begin{eqnarray*}
N_{J_\alpha}(X^\#_p,Y_p^\#) &=& [-\alpha^{-1}(X),-\alpha^{-1}(Y)]_p -
[X^\#,Y^\#]_p \\ & & \quad -                    J_\alpha \, [X^\#,  - \alpha^{-1}(Y)]_p - J_\alpha \, 
[-\alpha^{-1}(X),Y^\#]_p \ ,
\end{eqnarray*}
where by slight abuse of notation we denote by $\alpha^{-1}(X)$ the vector field $p \mapsto \alpha^{-1}_p(X)$. 
\\

We examine each term individually:
\begin{enumerate}[label=(\roman*)]
        \item We have $[X^\#, Y^\#]_p = ([X,Y])_p^\#  \  \in V_p$, so this term is vertical. \\

        \item We recall that the Lie bracket $[\zeta,\xi]_p$ of two vector fields may be expressed in terms of the one-parameter flow $\phi^t_\zeta$ of $\zeta$ by the equation
\[
	[\zeta,\xi]_p = \dtnull (\phi^t_\zeta)_{*,p}^{-1} ( \xi_{\phi^t_\zeta(p)} ) \ . 
\]
The one-parameter flow of $X^\#$ is given by $\phi^t_{X^\#} (p) = p e^{tX} = R_{e^{tX}}(p)$. 
With the required equivariance of $\alpha$ this gives:
\begin{eqnarray*}
[X^\#,\alpha^{-1}(Y)]_p & = & \dtnull (R_{e^{tX}})_{*,p}^{-1} \ (\alpha^{-1}_{p
e^{tX}} (Y) ) \\
                     & = & \dtnull (R_{e^{-tX}})_{*,pe^{tX}} \ (\alpha^{-1}_{p
e^{tX}} (Y) ) \\
                & = & \dtnull (\alpha_p^{-1} \circ \ad_{e^{tX}}) (Y) \\
                & = & \alpha_p^{-1} (\dtnull \ad_{e^{tX}} (Y)) \\
                & = & \alpha_p^{-1} ([X,Y]) \ .
\end{eqnarray*}
In particular, this vector is horizontal and thus
\begin{eqnarray*}
- J_\alpha [X^\#,  - \alpha^{-1}(Y)]_p & = & (^\# \circ \, \alpha)( [X^\#,  -
\alpha^{-1}(Y)]_p ) \\
                        & = & (^\# \circ \,  \alpha \circ \alpha^{-1}_p )([X,Y]) \\
                        & = & ([X,Y])^\#_p \ .
\end{eqnarray*}
Analogously
\begin{eqnarray*}
- J_\alpha [ -\alpha^{-1}(X), Y^\#]_p & = & ([X,Y])^\#_p \ .
\end{eqnarray*}
These terms are vertical. \\

\item 
The vector fields $\alpha^{-1}(X)$, $\alpha^{-1}(Y)$ are horizontal. Therefore we get
\begin{align*}
F_A(\alpha^{-1}(X), \alpha^{-1}(Y)) & = d\omega_A(\alpha^{-1}(X), \alpha^{-1}(Y)) \\
	& = -\omega_A((\alpha^{-1}(X), \alpha^{-1}(Y)) \ . 
\end{align*}
This yields the vertical component of the Lie bracket of $\alpha^{-1}(X)$ and  $\alpha^{-1}(Y)$:
\begin{eqnarray*}
        \pi_V ([\alpha^{-1}(X), \alpha^{-1}(Y)]_p) & = - ( F_A
(\alpha^{-1}(X), \alpha^{-1}(Y)))^\#_p \ . 
\end{eqnarray*}\\
\end{enumerate}
Using these computations, we can express the horizontal and vertical components of the vector field $N_{J_\alpha}(X^\#,Y^\#)$:
\begin{equation}\label{nj1}
\begin{split}
N_{J_\alpha}(X^\#_p,Y_p^\#) =&
\underbrace{\pi_A [\alpha^{-1}(X),\alpha^{-1}(Y)]_p }_{\in A}
            \\  & + \underbrace{\left([X,Y] - 
           F(\alpha_p^{-1}(X),\alpha_p^{-1}(Y)) \right)^\#_p}_{\in V}.
\end{split}
\end{equation}
We claim that 
\[
	\alpha ( \pi_A [\alpha^{-1}(X),\alpha^{-1}(Y)] ) = - (d_A \alpha)\, (\alpha^{-1}(X),\alpha^{-1}(Y)) \ .
\]
In fact, 
\begin{eqnarray*}
(d_A \alpha)(\alpha_p^{-1}(X),\alpha_p^{-1}(Y)) & = & d \alpha (\alpha_p^{-1}(X),\alpha_p^{-1}(Y)) \\	
& = & 
        \alpha_p^{-1}(X). \underbrace{\alpha(\alpha^{-1}(Y))}_{=Y}
        -\alpha_p^{-1}(Y). \underbrace{\alpha(\alpha^{-1}(X))}_{=X}
        \\
	& & - \alpha( [\alpha^{-1}(X),\alpha^{-1}(Y)]_p) \\
        & = & -\alpha( [\alpha^{-1}(X),\alpha^{-1}(Y)]_p) \\
	& = & -\alpha( \pi_A [\alpha^{-1}(X),\alpha^{-1}(Y)]_p) \ ,
\end{eqnarray*}
where we have made use of the fact that $\alpha$ vanishes on vertical vectors. 
\\

Recalling that $\alpha$ is a pointwise isomorphism $A_p \to \k$, we conclude that the vertical component of $N_{\Ja}(X^\#, Y^\#)$ vanishes for all $X,Y \in \k$ if and only if $d_A \alpha = 0$, noticing that the elements of the form $\alpha_p^{-1}(X)$, with $X \in \k$, generate the horizontal tangent space in $p$. 
\\

The vertical component of $N_{\Ja}(X^\#, Y^\#)$ vanishes if and only if 
\begin{eqnarray}\label{cond1}
 F_A\ (\alpha_p^{-1}(X),\alpha_p^{-1}(Y)) = [X,Y]
\end{eqnarray}
for all $X,Y \in \k$, which in the introduced hybrid notation is equivalent to 
$$
F_A = \frac{1}{2}\ [\alpha \wedge \alpha] \ .
$$
This achieves the proof of the Theorem \ref{main_prop}. \qed

\section{Integrability and induced metrics on the base}
The integrability equations are presumably difficult to solve in general. However we will get solutions in geometrical terms by the following idea: We endow the base $M$ with a metric $g_\alpha$ induced by $\alpha$  and an $\ad$-invariant inner product on the Lie algebra $\k$. We then get interesting correspondances between the Riemannian geometry of $(M,g_\alpha)$ and the integrability conditions. In particular, we obtain necessary conditions in general, including non-positive sectional curvature, and we get a complete solution to the integrability conditions for $K=SU(2)$ or $K=SO(3)$ if $M$ is a hyperbolic 3-manifold.
\\

We shall first review some notation. We know \cite{KN1} that there is an isomorphism between the $\ad$-invariant tensorial forms on $P$ with value in $\k$ and the forms on $M$ with value in the associated bundle $\ad(P)$. We denote this isomorphism by
$$
   ^\vee : A^k_{\ad}(P,\k) \nach A^k(M, P\times_{\ad} \k)
$$
which is given by the formula
$ ^\vee(\beta)_{\pi(p)}(u_1, \dots, u_k) = [p;\beta_p (\tildeu_1, \dots, \tildeu_k))]$, where the $\tildeu_i$ are horizontal liftings of the vectors $u_i$ and where $[p;X]$ denotes the class of $(p,X) \in P\times \k$ in the adjoint bundle.
We denote $^\wedge$ the inverse of $^\vee$. For $\gamma \in  A^k(M, P\times_{\ad} \k)$ we have the formula $(\w \gamma)_p(\lambda_1, \dots, \lambda_k) = \bar{\gamma}(\pi_* \lambda_1, \dots, \pi_* \lambda_k)$ such that $[p; \bar{\gamma}(\pi_* \lambda_1, \dots, \pi_* \lambda_k)] = \gamma(\pi_* \lambda_1, \dots, \pi_* \lambda_k)$. 

As $K$ was supposed to be a compact Lie group, there exists an $\ad$-invariant metric, that means an inner product $\left<.,.\right>$ satisfying
$$
\left<\ad_g X , \ad_g Y \right>\  =\  \left< X,Y\right>
$$
for all $g \in K$ and $X,Y  \in \k$. By invariance this descends to the associated bundle $\ad(P) = P \times_{\ad} \k$ which thereby naturally obtains the structure of an Euclidean vector bundle. We shall denote by $\t$ the corresponding metric. By assumption $^\vee \alpha$ is an isomorphism $TM \nach \ad(P)$. We can pull back the metric $\t$ to $TM$, obtaining a Riemannian metric 
$\ga$ on $M$. $\val$ becomes than an isometry of Euclidean vector bundles. 
\\

The connection $A$ on $P$ induces a linear connection $\nabla^A$ on $\ad(P)$. In what follows we use two of the equivalent notions of $\nabla^A$ \cite{KN1} which we briefly recall. Let $\sigma$ be a section of $\ad(P)$ (this is an element of $A^0(M,\ad(P))$). Then
$$
 \nabla^A \sigma = (\v \ \circ \ d_A \ \circ \ ^\wedge)(\sigma) \ .
$$
Another equivalent definition is constructed as follows: Let $\tau=(x(t))$ be a curve on $M$. The point $\tau^{t+h}_t (p)$ is then the end point of the unique horizontal lifting of $\tau$ beginning at $p \in \piinv(x(t))$ and ending at $\piinv(x(t+h))$. This gives the notion of parallel displacement on $P$. This construction carries over to parallel displacement in the associated bundle by the formula
$$
\tau^{t+h}_{t} ([p;X]) := [\tau_t^{t+h}(p); X]
$$
(which is defined independtly of representatives). For a section $\sigma$ of the associated bundle one then has
$$
 \nabla^A_{\dot{x}_t} \sigma = 
\left. \frac{d}{dh} \right|_{h=0} \underbrace{\tau^{t}_{t+h} (\sigma(x(t+h)))}_{\in \ad(P)_{x(t)} \ \forall h} .
$$

Given the connection $\nabla^A$ on $\ad(P)$, there is a natural operator $d^{\nabla^A}$ on forms with value in the associated bundle, generalising the usual exterior derivative. It is defined by
\begin{eqnarray*}
d^{\nabla^{A}} (\eta \tensor s) &=& d\eta \tensor s + (-1)^{\deg(\eta)} \eta \wedge d^{\nabla^{A}} s \\
        d^{\nabla^{A}} s & = & \nabla^A s
\end{eqnarray*}
for $\eta \in A^k(M,\R)$ and $ s \in A^0(M, \ad(P))$.
\\

Let  $\mathfrak{ad} : \k \rightarrow \gl(\k)$ be the Lie algebra homomorphism corresponding to the adjoint represenation $\ad: K \rightarrow \Gl(\k)$, i.e. $\mathfrak{ad}$ is the derivativative of $\ad$ at the identity and is therefore just the endomorphism induced by the Lie bracket. 
Let $\eta \in A_{\ad}^1(P,\k)$ be a tensorial 1-form of type $\ad$. Then one has
\begin{eqnarray}\label{d_A}
        d_A \eta = d \eta + \mathfrak{ad}(\omA) \wedge \eta \, ,
\end{eqnarray}
or, more explicitely, for two vectors $\zeta, \xi$ 
\begin{eqnarray*}
        d_A \eta (\zeta,\xi)   &=& 
                d \eta(\zeta,\xi) + \mathfrak{ad}(\omA(\zeta))  \eta(\xi)
                                -   \mathfrak{ad}(\omA(\xi))  \eta(\zeta) \\
        & = &   d \eta(\zeta,\xi) + [\omA(\zeta),  \eta(\xi)] 
                        - [\omA(\xi),  \eta(\zeta)] \ .
\end{eqnarray*}
The curvature of the connection $\nabla^A$, applied to a section $\sigma$ of
 $\ad(P)$, is given by
\begin{eqnarray}\label{FnablaA}
F^{\nabla^A} (\sigma) & := & \dnablaA \circ \dnablaA (\sigma) \nonumber \\ 
       & =& \v \circ d_A \circ d_A \circ \w (\sigma) \nonumber \\
       & = & \v \circ \mathfrak{ad}(F_A) \circ \w (\sigma) \ ,
\end{eqnarray}
using (\ref{d_A}).

With these notions, we can naturally define a linear connection $\nabla$ on the tangent bundle $TM$ by pulling back the connection $\nabla^A$ on $\ad(P)$ via $\val$: 
$$
\nabla := (\val)^{-1} \nabla^A :=  (\val)^{-1} \circ \nabla^A \circ (\val \ ) .
$$
As $\ad$ operates by isometries on $\k$ the connection $\nabla^A$ is compatible with the metric $\t$ . As $\val$ is an isometry, $\nabla$ is compatible with the metric $g_\alpha$. We recall that the Levi-Civita connection on the tangent bundle of a Riemannian manifold is the unique connection that is compatible with the metric and that has vanishing torsion. 
\begin{prop}\label{torsion}
  The torsion $T_\nabla$ of the connection $\nabla$ constructed above vanishes if and only if $d_A \alpha = 0$. If this is the case, then $\nabla$ is the Levi-Civita-connection of $(M,g_\alpha)$. 
\end{prop}
\proof
Fix $p \in P$ and let $x:=\pi(p)$, $u_x \in T_xM, \ v_x \in T_xM$. Take $u,v$ arbitrary vector fields that are extensions of $u_x,v_x$, and take $\tildeu, \tildev$ horizontal liftings. Note that $\w((\val)(v))(p) = \alpha_p(\tildev)$.
Thus one has 

\begin{eqnarray*}
(\val)(T_{\nabla} (u_{x},v_{x})) & = & 
\nablaA_{u_x} (\val)(v) - \nablaA_{v_x} (\val)(u) - (\val) ( [u,v]_x ) \\
& = & [(\cdAh)((\val)(v))](u_x) - [(\cdAh)((\val)(u))](v_x) \\ & & - (\val) ( [u,v]_x ) \\
& = & \v (d_A \alpha(\tildev)) (u_x) - \v (d_A \alpha(\tildeu)) (v_x) -  (\val) ( [u,v]_x ) \\
& = & [p; d_A (\alpha(\tildev)) (\tildeu_p) -  d_A (\alpha(\tildeu)) (\tildev_p) - \alpha(\widetilde{[u,v]}_p) ] \, . \\
\end{eqnarray*}
As we easily see that $\pi_* \circ [\tildeu,\tildev] = \pi_* \circ \widetilde{[u,v]}$, and as $\alpha$ is a tensorial form, we get
$\alpha(\widetilde{[u,v]}_p) = \alpha( [\tildeu,\tildev]_p)$. 
This yields finally that 
\begin{eqnarray*}
(\val)(T_{\nabla} (u_{x},v_{x})) & = & [p; (d_A \alpha)(\tildeu_p,\tildev_p)] .
\end{eqnarray*}
The claim follows. \qed

\begin{prop} Suppose the above integrability equations are satisfied. Then the Riemannian curvature $R_\nabla$ of $\nabla$ is given by the formula
$$
R_\nabla (u_x,v_x) w_x = (\val)^{-1} \left(\left[p;[[\alpha(\tildeu_p),\alpha(\tildev_p)], \alpha(\tildew_p)]\right] \right)
$$
with the above notations. Its sectional curvature is always non-positive. It is strictly negative if the dimension of the maximal torus of $K$ is 1.
\end{prop}
\proof
The curvature formula follows from the above formula (\ref{FnablaA}) and the definition of $R_\nabla$ as $d^\nabla \circ d^\nabla$, as well as the explicit form of $F_A$ given by the integrability equations. 
For the sectional curvature we notice first that:
\begin{eqnarray}\label{formule_t}
\left<[Z,X],Y\right> + \left<X,[Z,Y]\right> = 0
\end{eqnarray}
for all $X,Y,Z \in \k$. 
This follows from the $\ad$-invariance by derivation 
\begin{eqnarray*}
0 & = & \dtnull \left< \ad_{e^{tZ}}X, \ad_{e^{tZ}}Y\right> \\
  & = & \left< \dtnull \ad_{e^{tZ}}X,Y\right> + \left<X, \dtnull \ad_{e^{tZ}}Y\right> \\
  & = & \left<[Z,X],Y\right> + \left<X,[Z,Y]\right> .
\end{eqnarray*}

Now the sectional curvature $K(\langle \, u,v \rangle )$ of the plane $ \langle u,v \rangle$  spanned by the orthonormal vectors $u,v \in T_xM$ is given by
$$
K(\langle \, u,v \rangle ) = g_\alpha(u, R_\nabla (u,v)v)
$$
By the curvature formula, the above formula (\ref{formule_t}) and the definition of $g_\alpha$ 
we get:
\begin{eqnarray*}
K(\langle \, u,v \rangle ) & = & 
\t\left((\val) (u),[p; [[\alpha(\tildeu_p),\alpha(\tildev_p)], \alpha(\tildev_p)]]) \right) \\
        & = & \left< \alpha(\tildeu_p), [[\alpha(\tildeu_p),\alpha(\tildev_p)], \alpha(\tildev_p)] \right> \\
        & = & - \left< [\alpha(\tildeu_p), \alpha(\tildev_p)] , [\alpha(\tildeu_p), \alpha(\tildev_p)] \right> \\
        & \leq & 0 \ . 
\end{eqnarray*}
\qed
\subsection{Special cases $K=SU(2)$ and $K=SO(3)$}
We shall show here that the integrability conditions are related to hyperbolic geometry.
Recall that the Lie groups $SU(2)$ and $SO(3)$ are locally isomorphic, hence they have isomorphic Lie algebras. There is an $\ad$-invariant metric on $\su(2) \cong \so(3)$ given by
$$
\left<X,Y\right> := \tr (XY^*) = - \tr (XY) \quad \forall X,Y \in \su(2) \ .
$$
\begin{theorem}
Let $P \stackrel{\pi}{\nach}M$ be an $SU(2)$- or $SO(3)$- principal fibre bundle and $\val : TM \nach \ad(P)$ an isometry, where $\ad(P)$ possesses the metric induced by the above metric on $\su(2)$, respectively $\so(3)$. If the almost-complex structure $J_\alpha$ is integrable, then $(M,g_\alpha)$ has constant negative sectional curvature.
\end{theorem}
\proof
We recall the Pauli-matrices
\begin{eqnarray*}
\sigma_1=\left(\begin{array}{cc} 0 & 1 \\ 1 & 0 \end{array} \right) \hspace{1 cm}
\sigma_2=\left(\begin{array}{cc} 0 & -i \\ i & 0 \end{array} \right)\hspace{1 cm}
\sigma_3=\left(\begin{array}{cc} 1 & 0 \\ 0 & -1 \end{array} \right)
\end{eqnarray*}
which satisfy
$$
        \sigma_i \sigma_j = \delta_{ij} 1 + i \epsilon_{ijk} \sigma_k \ . 
$$
Consequently, the matrices 
$$
X_i := - \frac{i \sigma_i}{\sqrt{2}}
$$
form an orthonormal basis of $(\su(2),\left<.,.\right>)$ and satisfy the relation
\begin{eqnarray*}
[X_i,X_j] = \sqrt{2} \epsilon_{ijk} X_k \ .
\end{eqnarray*}
If now the vectors $u,v \in T_x M$ form an orthonormal basis of the 2-dimensional plane $\langle \, u,v \rangle  \  \subseteq T_xM$, then the elements
\begin{eqnarray*}
        [p;Z]:= (\val)(u) \ \ \ [p;W]:= (\val)(v)
\end{eqnarray*}
are orthonormal in $\ad(P)_x$ and thus $Z,W$ are with respect to $(\su(2),\left<.,.\right>)$
The decompositions with respect to the above orthonormal basis  
$W =: w^i X_i$, $Z=: z^i X_i$ implies
\begin{eqnarray*}
        \left<W,W\right> & = \sum_{i=1}^{3} (w^i)^2 &= 1 \\
        \left<Z,Z\right> & = \sum_{i=1}^{3} (z^i)^2 &= 1 \\
        \left<Z,W\right> &= \sum_{i=1}^{3} z^i w^i & = 0 \ .
\end{eqnarray*}
We use the formula for the sectional curvature in the proof of the last Theorem now:
\begin{eqnarray*}
K(\langle \, u,v \rangle ) & = & - \left<[W,Z],[W,Z]\right> \\
& = & - 2 w^i z^j w^k z^l \ \left<\epsilon_{ijm} X_m, \epsilon_{kln} X_n \right> \\
& = & -2  w^i z^j w^k z^l \ \epsilon_{ijm}\epsilon_{kln} \delta^{mn} \\
& = & -2  w^i z^j w^k z^l \ (\delta_{ik} \delta_{jl} - \delta_{il} \delta_{jk}) \\
& = & -2 \left<W,W\right> \left<Z,Z\right> + 2 \left<W,Z\right> \\
& = & -2
\end{eqnarray*}
\qed
This result leads to the question whether on all hyperbolic 3-manifold there is a principal fibre bundle with an integrable almost complex structure. This is indeed the case: 
\begin{theorem}\label{solution hyperbolic case}
  Let $(M,g)$ be a hyperbolic oriented 3-manifold. Then the $SO(3)$-principal fibre bundle $P_{SO(M)}$ of orthonormal oriented frames has a natural integrable almost-complex structure.
\end{theorem}
As a consequence every such principal fibre bundle is naturally a complex 3-dimensional manifold. \\
\\
\proof
We begin by defining $\val$ which becomes an isometry $TM \stackrel{\cong}{\nach} \ad(P)$. To this end we first define 
$$
L : \R^3 \nach \so(3)
$$
where the image $L_x$ of an element $x \in \R^3$ is the endomorphism defined by $L_x(y) := x \times y$ with $\times$ the usual vector product. This is indeed well defined since $L_x$ is anti-symmetric and traceless. $L$ is a Lie algebra isomorphism and has the important equivariance property given by the following commutative diagram
\begin{eqnarray}\label{L for so(3)}
\begin{diagram}
\node{SO(3) \times \R^3 } \arrow[2]{e,t}{can. rep.}
\arrow[1]{s,l}{id\times L}  \node[2]{\R^3}\arrow[1]{s,r}{L}  \\
\node{SO(3) \times \so(3)} \arrow[2] {e,b}{adj. rep.} \node[2]{\so(3)}
\end{diagram} 
\end{eqnarray}
If now we define the metric on  $\so(3)$ which is one half the usual metric given by the trace formula, then $L$ becomes an isometry.\\

We recall that the tangent bundle $TM$ is isometrical to the associated bundle of the orthonormal oriented principal fibre bundle $P_{SO(M)}$ with respect to the canonical representation of $SO(3)$ which leaves the usual metric on $\R^3$ invariant.
$$
TM \cong P_{SO(M)} \times_{can} \R^3
$$
We identify these euclidean vector bundles.
Thanks to the above commutative diagram we may safely define:
\begin{eqnarray*}
  \val : \ TM = \Psom \times_{can} \R^3 \nach \Psom \times_{\ad} \so(3) = \ad(P)
\end{eqnarray*}
by $\val([p,y]_{can}) = [p,L(y)]_{\ad}$.
Let $\nabla$ denote the Levi-Civita-Connection of $(M,g)$, and let $A$ be the associated connection in $\Psom$. Let $\nabla^A$ denote the induced connection on $\ad(P)$. 
\begin{lemma}
$$  
\nabla^A = (\val) (\nabla) := (\val) \circ \nabla \circ (\val)^{-1} 
$$
\end{lemma}
\qed
This Lemma is true as one finds by going to local sections and by the above diagram for $L$.\\

The Lemma implies that we have $(\val)^{-1}(\nabla^A) = \nabla$. By Proposition \ref{torsion} and the fact that the Levi-Civita connection $\nabla$ is torsion free, we get $d_A \alpha = 0$, yielding the first of the integrability equations.\\

The preceeding Lemma implies $F^{\nabla^A} = F^{(\val) (\nabla)}$. As the connection $F^\nabla$ is equal to the Riemannian curvature $R_\nabla$, we get the relation
$$
F^{\nabla^A} = ( \val ) \circ R_\nabla \circ (\val)^{-1} \ . 
$$
By assumption the sectional curvature is equal to a constant, say $\kappa$, and this yields the formula
$$
R_\nabla (v,w) u = \kappa \ (g(u,w)v - g(u,v)w) \ . 
$$
for $u,v,w \in T_xM, \ x \in M$. Using this, we get the following formula for the curvature $F^{\nabla^A}$, used further down. Let $s \in \ad(P)_x$ with $s= (\val)(u)$.
Then
\begin{eqnarray*}
F^{\nabla^A}(v,w) s & = & (\val)(R_\nabla(v,w) u ) \\
& = & (\val)(\kappa \  (g(u,w)\ v - g(u,v) w)) \\
& = & \kappa \ \left( \t(s,(\val)(w)) \ (\val)(v) - \t(s,(\val)(v)) \ (\val)(w)   \right) \ . 
\end{eqnarray*}
If now we define $U, V, W \in \so(3)$ by the formulae $s=:[p;U],\  (\val)(w) = [p;W]$, and $ (\val)(v) = [p;V]$, we get the final formula
\begin{eqnarray}\label{FnA}
F^{\nabla^A}(u,w) \ s = \kappa \ [p; \left<U,W\right> V - \left<U,V\right> W ] \ .
\end{eqnarray}

Our aim is to establish the formula
\begin{eqnarray}\label{tofind}
  F^{\nabla^A} = \v \circ \, \mathfrak{ad} (\frac{1}{2} [\alpha \wedge \alpha]) \, \circ \, \w \,
\end{eqnarray}
as the Theorem will follow, because this will imply that $F_A = \frac{1}{2} [\alpha \wedge \alpha]$, yielding the second integrability equation. That this indeed implies $F_A = \frac{1}{2} [\alpha \wedge \alpha]$ follows from the formula (\ref{FnablaA}) and the fact that $\mathfrak{ad}$ is injective, since surely $\so(3)$ has trivial center. We calculate now the right hand side of the last equation, using the above notations:
\begin{eqnarray*}
 (\v \, \circ \, \mathfrak{ad} (\frac{1}{2} [\alpha \wedge \alpha]) \, \circ \,  \w)(v,w) s 
 & = & (\v \circ \mathfrak{ad}(\frac{1}{2} [\alpha \wedge \alpha]))(v,w) [p;U] \\
 & = & [p; \mathfrak{ad}(\frac{1}{2} [\alpha \wedge \alpha])(\tildev, \tildew) U] \\
 & = & [p; [[V,W],U]] \ . 
\end{eqnarray*}
Comparing this with (\ref{FnA}), we get the equation in question (\ref{tofind}) if and only if we have
\begin{eqnarray*}
  [[V,W],U] = \kappa ( \left<U,W\right> V - \left<U,V\right> W )
\end{eqnarray*}
for all $U,V,W \in \k$. This is really the case for $\kappa=-1$ since then this equation reduces to the well-known vector-product identiy in $\R^3$ (via the Lie algebra isomorphism $L$)  
\begin{eqnarray*}
  (\mathfrak{v}\times\mathfrak{w})\times\mathfrak{u} 
= (\mathfrak{u}\cdot\mathfrak{v})\mathfrak{w} 
- (\mathfrak{u}\cdot\mathfrak{w})\mathfrak{v} \, 
\end{eqnarray*}
with $v=[p,\mathfrak{v}]$, $p$ an orthonormal frame and $\mathfrak{v} \in \R^3$, and likewise $w =[p,\mathfrak{w}]$, $u=[p,\mathfrak{u}]$. Of course, we may always rescale the metric so that we have $\kappa = -1$. 
\\
This completes the proof. \qed

\subsection{Further examples of the required topological setting}
In the previous construction the existence of a 1--form $\alpha$, defining an isomorphism 
$TM  \to \ad(P)$, where $P \to M$ was the principal oriented orthonormal frame bundle of $M$, can be seen as an `accident of low dimension'. Indeed, the dimension of the Lie group $SO(n)$ is strictly larger than $n$ for all $n \geq 4$. However, we may again be in a similar setting if the structure group of the tangent bundle $TM$ admits a reduction from $SO(n)$ to a subgroup $K$ which has dimension equal to the dimension of the base manifold $M$. 

More concretely, suppose we have a 6--dimensional Riemannian manifold $M$ that has a reduction of its structure group to $K=SO(3) \times SO(3)$. Let $P \to M$ be the $K$--principal bundle obtained from the principal oriented orthonormal frame bundle $P_{SO(M)} \to M$ by reduction to the structure group $K$. The construction (\ref{L for so(3)}) in the proof of Theorem \ref{solution hyperbolic case} extends to the situation here yielding an isomorphism
\[ TM \to P_{\ad_+ \times \ad_-} \so(3)_+ \oplus \so(3)_- 
\ ,\]
thereby giving rise to an isomorphism $\alpha$ as required for the construction of the almost complex structure $\Ja$. By Berger's classification \cite{B} of simply connected Riemannian manifolds with special holonomy, we may not expect to have a reduction of the structure group to $SO(3) \times SO(3)$ for a simply connected 6--manifold $M$ that is not a symmetric space or a locally product space. However, for locally product spaces there are many examples that fall into this class, for instance, all closed 6--manifolds that fibre over a 3--manifold. Presumably, this class includes cases where one of the $SO(3)$--bundles $\ad_\pm(P)$ is non-trivial. 

Another large class of examples where our construction gives an almost complex structure can be obtained from parallelisable manifolds $M$. In fact, if $n$ is the dimension of $M$, then any $T^n$--principal bundle $P \to M$ is such that the bundles $TM$ and $\ad(P)$ are isomorphic, and so one may obtain a class $\alpha$ as required in our construction. There are of course many situations where $P$ is a non-trivial bundle. In fact, for any collection of complex line bundles $L_i \to M, \, i = 1, \dots, n$, the principal complex frame bundle of $L_1 \oplus \dots \oplus L_n$ with its natural reduction of the structure group from $U(n)$ to $S^1 \times \dots \times S^1 = T^n$ is such an example. These are in 1--to--1 correspondance with classes $e_1, \dots ,e_n \in H^2(M;\Z)$ via the corresponding first Chern classes.

One may expect more examples related to `special holonomy'. 

\section{$J_\alpha$-holomorphic vector fields and completeness}
In this section we stay in the setting of a principal bundle $P \to M$ with compact structure group $K$ endowed with a connection $A$ on $P$, and with an almost complex structure $J_\alpha$ determined by $A$ and the 1--form $\alpha$. We assume an $\ad$--invariant inner product on the Lie algebra of $K$ turning $\ad(P)$ into an Euclidean vector bundle, and we suppose $M$ is endowed with a Riemannian metric $g_\alpha$ turning the induced map $^\vee \alpha: TM \to \ad(P)$ into an isometry. We also suppose having the unique connection $\nabla$ on $TM$ such that $(^\vee \alpha)(\nabla) = \nabla^A$. 
We will not suppose $J_\alpha$ to be integrable in this section, but only that $d_A \alpha = 0$. By Proposition \ref{torsion} above, this implies that the connection $\nabla$ is the Levi-Civita connection on $TM$. 

We will show here that there is an equivalence between the complete integrability of the naturally occuring $\Ja$-holomorphic vector fields $X^\# - i \ J_\alpha X^\#$ and the geodesic completeness of the Riemannian manifold $(M,g_\alpha)$. 
Here the vector field $X^\# - i \ J_\alpha X^\#$ is said to be completely integrable if for any $p \in P$ there is a map $\phi : \C \nach P$ which is $\Ja$-holomorphic, i.e. such that
$$
\phi_* \left(\frac{\partial}{\partial z} \right) = \frac{1}{2} (X^\# - i \ J_\alpha X^\#)
 \quad (X \in \k)
$$
for $\frac{\partial}{\partial z} = \frac{1}{2}(\frac{\partial}{\partial x } - i \frac{\partial}{\partial y })$, 
and which satisfies $\phi(0) = p$.

\begin{prop}
The $J_\alpha$-holomorphic vector fields (sections of $T^{1,0}P$)
$$
X^\# - i \ J_\alpha X^\# \quad 
$$
are completely integrable for any $X \in \k$ if and only if $(M,g_\alpha)$ is geodesically complete.
\end{prop}
\proof
The correspondance is based on the following lemma:
\begin{lemma}
1. Let $X \in \k$ and denote $\phi^t_{\Ja X^\#}$ the flow of the vector field $\Ja X^\#$, defined for all $p \in P$ at least in a neighbourhood of $0 \in \R$. Then for all $p \in P$ the curve $c(t):= \pi \circ \phitJaX (p)$ is a geodesic for the values of $t$ where this is defined.\\
2. Let $c(t)$ be a geodesic on M defined at least in a neighbourhood of $0 \in \R$ . Then there is an element $X \in \k$ such that
 $c(t):= \pi \circ \phitJaX (p)$ for all $p \in \piinv(c(0))$ for these values of $t$.
\end{lemma}
\proof
1. We shall abbreviate our notation by $\phi^t(p) := \phitJaX (p)$.
%$$
%\dot{c}(t) = \pi_{*,\phi^t(p)} \JaX _{\phi^t(p)}
%$$
One then gets by the definition of $c(t)$ 
\begin{eqnarray*}
  \nabla_{\cdt} \cdt & = & (\val)^{-1} (\nabla^A_{\cdt} (\val)(\cdt)) \\
  & = & (\val)^{-1} (\nabla^A_{\cdt} (\val)(\pi_{*,\phi^t(p)} \JaX _{\phi^t(p)}))
  \\
  & = & (\val)^{-1} (\nabla^A_{\cdt} [\phi^t(p); -X] \\
  & = & (\val)^{-1} \left. \frac{d}{dh} \right|_{h=0} 
  \underbrace{c^t_{t+h} [\phi^{t+h}(p);-X]}_{=[\phi^t(p);-X]} \\
  & = & 0
\end{eqnarray*}
since in fact $\phi^t(p)$ is already a horizontal lifting of $c(t)$. By definition, this means that the curve $t \mapsto c(t)$ is a geodesic. \\

2. Let $t \mapsto c(t)$ be a geodesic on $M$, defined in a connected neighbourhood of $0$. Let $p \in \piinv(c(0))$. Then there is an
$ X \in \k$ such that
\begin{eqnarray}\label{condinitial}
  \pi_{*,p} \JaX = \dot{c}(0) \ .
\end{eqnarray}
Set $\psi(t) := \pi \circ \phitJaX (p)$. As we have seen in the proof of (1.) we get
$$
\nabla_{\dot{\psi}(t)} \dot{\psi}(t) = 0
$$
as well as $\nabla_{\cdt} \cdt = 0$. By the equation (\ref{condinitial}) they have the same initial condition. Now by the uniqueness of the solution of a Cauchy problem we get the claim, at least for $t$ in a certain neighbourhood of $0$. \qed \\
We continue with the proof of the proposition.
Let $t \mapsto c(t)$ be a geodesic on $M$, defined in a neighbourhood of $0$. Let $X \in \k$ be as in the formula (\ref{condinitial}). The vector field $X^\# - i \ J_\alpha X^\#$ is completely integrable by assumption. This means that there is a map $\phi : \C \nach P$ which is $\Ja$-holomorphic, therefore satisfying
$$
\phi_* \left(\frac{\partial}{\partial z} \right) = \frac{1}{2} (X^\# - i \ J_\alpha X^\#)
 \quad (X \in \k)
$$
for $\frac{\partial}{\partial z} = \frac{1}{2}(\frac{\partial}{\partial x } - i \frac{\partial}{\partial y })$, and which satisfies $\phi(0) \in \piinv(c(0))$. The map $y \mapsto \phi(0+iy)$ is an integral curve of the horizontal vector field $\JaX$. This curve, composed with the projection $\pi$ to $M$, is the geodesic $c$, as we have seen before. This is true at least in the neigbourhood where $c$ is defined. But the map $ y \mapsto \phi(0+iy)$ has a continuation to the whole of $\R$, so $c$ posseses a continuation to $\R$ as a geodesic.
\par

Let now $M$ be geodesically complete. Consider the map 
$$
s+it \mapsto \phitJaX \circ \phi^s_{X^\#} (p)
$$
for an arbitrary but fixed $p \in P$. Since the Lie bracket of the two vector fields $X^\#$ and $\JaX$ vanishes (this follows from the computation of the second term of Nijenhuis tensor in the proof of Theorem \ref{main_prop} above), we know that the two associated flows commute, and thus the above map is $\Ja$-holomorphic. This implies for an allowed value $t$ that we can vary in the $s$-direction on the whole real line. Thus the maximal domain of definition is always of the form $I \times \R$ with $I$ open and connected in $\R$, containing $0$. We actually want to show that $I = \R$. We denote $\phi(t):= \phitJaX(p)$ and know that $t \mapsto c(t):= (\pi \circ \phi)(t)$ is a geodesic on $M$, defined for $t \in I$. Suppose $t_0 := \sup (I) < \infty$. By the geodesic completeness of $M$ we get a geodesic 
$$
\bar{c} : \R \nach M
$$
with $\bar{c}|_I = c$. Let $\beta$ be a horizontal lifting of $\bar{c}$ such that $\beta(0) = \phi(0)$ (the lifting $\beta(t)$ exists for all $t$ in the domain of definition of $\bar{c}$, see \cite{KN1} p. 69.). We then have by construction $\beta|_I = \phi$, and in particular
$$
\dot{\beta}(t) = \JaX_{\beta(t)} 
$$
for $t<t_0$. But as all maps are continous, we get by passing to the limit 
$t \rightarrow t_0, t<t_0$ that 
$$
\dot{\beta}(t_0) = \JaX_{\beta(t_0)} \ .
$$
Thus we have $t_0 \in I$  and by maximality we get a contradiction. In the same manner we get $\inf(I) = -\infty$ and thus  morphic map extends to $\C$. Thus all vector-fields of the type $X^\# - i \JaX$ are completely integrable. \qed

\section{Transitive holomorphic right-action from the complexified group $K^\C$ and a converse of our construction}

We will suppose here the almost complex structure $\Ja$ is integrable and that the base manifold $(M,g_\alpha)$ is geodesically complete. Let $G = K^\C$ be the complexification of the compact Lie group $K$, see \cite{BtD}, \cite{Ho} for this notion. We will now show that there is a natural holomorphic action of $G$ on $P$ which is locally free and an extension of the action $\rho$ of the structure group $K$ on $P$. The new action is transitive provided $M$ is connected. Given this data, we finally obtain a natural $\g$-valued holomorphic 1-form on $P$, where $\g$ is the Lie algebra of $G$, the complexification of the Lie algebra $\k$. 
\\

We recall the KP-decomoposition of the complexification of a compact Lie group \cite{HH}: The holomorphic map $\g \to G$ given by $X+iY \mapsto \exp(X) \exp(iY)$ for $X,Y \in \k$, factoring through $K \times \k$, induces a diffeomorphism $K \times \k \to G$. It gives the unique decomposition $g=k\ e^{iX}$ of an element $g \in G$ with $k \in K$, and $ X \in \k$. This is just the usual polar decomposition for $G=Gl(n,\C)$. Via this, we now define a holomorphic map which is a priori no group action:
\begin{eqnarray*}
        \psi : \left\{ \begin{array}{ccc} P \times G & \longrightarrow & P \\
                 (p,k \ e^{iX}) & \longmapsto &  \phi^1_{\JaX} (pk) \ .
        \end{array} \right.
\end{eqnarray*}
That this is well defined in the case that the base manifold $(M,g_\alpha)$ is geodesically complete follows from the results of the last section.
Differentiating with respect to the second argument at the identity yields a section-valued map
\begin{eqnarray*}
  \psi_{*,e} : \g \nach \Gamma(TP) \ , \\
\end{eqnarray*}
which is the collection of all $(\psi_p)_{*,e}$. 
It satisfies
\begin{eqnarray*}
  \psi_{*,e}(X) & = & X^\# \\
  \psi_{*,e}(iX) & = & \JaX \ 
\end{eqnarray*}
for all $X \in \k$. 
Even more is true: The map $\psi_{*,e}$ is a homomorphism of (real) Lie-algebras (where the complexified Lie algebra $\g$ is, for the moment, considered as a real vector space) that interchanges the almost complex structure $\Ja$ and the canonical one on $\g$. In fact, that we have 
\[
	\psi_{*,e}[X,Y] = [\psi_{*,e} X, \psi_{*,e} Y] \ \ \text{and} \ \ 
	\psi_{*,e}[X,iY] = [\psi_{*,e} X, \psi_{*,e} iY]
\]
for all $X,Y \in \k$ follows from the computations in (i) and (ii) inside the proof of Theorem \ref{main_prop}, whereas the condition
\[
	\psi_{*,e}[iX,iY] = [\psi_{*,e} iX, \psi_{*,e} iY] 
\]
follows from the vanishing of the Nijenhuis tensor. 

The map $\psi_{*,e}$ becomes $\C$-linear after passing to the complexified tangent space of $TP$. By the second fundamental theorem of Lie there is thus an associated holomorphic {\em local} action $\phi$ of $G$ on $P$ with Lie-algebra homomorphism $\psi_{*,e}$, see for instance \cite{Ak}. By the associated one-parameter transformation groups one sees that $\phi$ has to coincide with $\psi$ where $\phi$ is defined. That $\psi$ defines a global action now follows from an identity theorem argument. We write $\psi = \rho^\C$ in what follows, indicating the extension of the already existing action $\rho$.
\\

Given the holomorphic action $\rho^\C: P \times G \to P$, we can define a holomorphic 1--form on $P$ with values in $\g$. In fact, partial differentiation of $\rho^\C$ in the second variable at the identity $g=e$ yields a biholomorphism, which in fact is a bundle isomorphism from the trivial bundle $\g \times P \to P$ to $TP$, 
\begin{equation*}
	\left. \frac{\partial \rho^\C}{\partial g}\right|_{g = e} : \g \times P \to TP \ .
\end{equation*}
We define $\omega \in \Omega^{1}(P,\g)$ to be the inverse of this map, followed by the projection $ \g \times P \to \g$. This is clearly a holomorphic 1--form, and its real part is given by the connection $1$--form, 
\[
\Re \, \omega = \omega_A \ .
\]

$G$ acts transitively on $P$ if $M$ is connected, as one sees by the fact that every two points on the base $M$ may then be joined by a geodesic. In general this geodesic is not unique and thus the action of $G$ is generally not free. More precisely, it is free if and only if the base manifold $(M,g_\alpha)$ is such that any two points are connected by a unique geodesic.
\\

We have now seen that in the case of geodesic completeness of the base $(M,g_\alpha)$, the complex structure on the manifold $P$ is in fact given as the quotient of a complex Lie group $G$ by a discrete subgroup $\Gamma$. This is in fact a more general phenomenon. By a result of Wang \cite{Wan}, any closed complex manifold with complex parallelisable tangent bundle is in fact given by the quotient of a complex Lie group by a discrete subgroup. However, the proof relies essentially on the compactness of the complex manifold, implying that certain complex-valued holomorphic functions on the manifold are necessarily constant by Liouville's theorem. 

Our result here may be seen as a generalisation of this result in the particular setting of an integrable almost complex structure $\Ja$ as studied here. Indeed, we don't have to require that $P$ is compact. In fact, we were only requiring the base manifold $(M,g_\alpha)$ to be geodesically complete which is a much weaker condition than requiring $M$ to be closed. 
\\

Conversely, one may ask how general our construction is in the situation that we have a complex parallelisable complex manifold $Q$ that is given as the quotient of a complex Lie group $K^\C$ by a discrete subgroup, and on which the Lie group $K$ acts freely. We will show now that in this situation the natural almost complex structure induced by the complex structure of the manifold $Q$ may be obtained as in our setting. 

For this we start with a `sample' construction. We assume that the complex Lie group $G$ is the complexification of the Lie group $K$. We may trivialise the tangent bundle $TG \to G$ as $\g \times G \to G$ by left-$G$-invariant vector fields defined by left $G$-translates of elements of the Lie algebra $\g = T_eG$. The group $G$ may be considered as a $K$-principal bundle $G \to G/K$ over the homogeneous space $G/K$. 

Any element of $\g = \k \tensor \C$ may be written as $X + iY$ with $X,Y \in \k$. We denote by $(X+iY)^\#$ the associated left-$G$-invariant vector field. We define $\alpha \in A^1_{\ad}(G;\k)$ by the following formula:
\[
	\alpha((X+iY)^\#) = -Y \ .
\]
With respect to $K$ this is really an $\ad$-equivariant tensorial one-form because we are considering $G$ to be a {\em right} $K$-principal bundle, with the natural action of $K$ on $G$ from the right. However, it is invariant under $G$ acting from the left. 
It is clear that the corresponding almost complex structure $\Ja$ coincides with the natural one on $G$, with complex charts coming, for instance, from the exponential map. 

We shall now assume that we have a transitive holomorphic action of $G$ on a connected complex manifold $Q$, with discrete stabiliser $\Gamma$ of some point in $Q$. In particular, $Q$ is biholomorphic to $\Gamma \backslash G$. We shall also assume that $K$ acts freely on $Q$. Because $\Gamma$ acts from the left on $G$, and because $\alpha$ is left-$G$-invariant, the form $\alpha$ descends to an $\ad$-equivariant tensorial one-form on $Q$, that we shall denote by $\bar{\alpha}$. It is clear now that the natural almost complex structure $J_{\bar{\alpha}}$ on $Q$, given as in our construction, coincides with the canonical almost complex structure on $Q$. We have thus given a sketch of proof of the following Proposition:

\begin{prop}
Suppose we have a connected complex manifold $Q$ with transitive holomorphic action of the complexification $G = K^\C$ of a compact Lie group with discrete stabiliser $\Gamma$ of some point in $Q$, such that the subgroup $K$ acts freely on $G$. Then the canonical almost complex structure on $Q$ is equal to the one given by our main construction $J_{\bar{\alpha}}$ induced from the sample situation on the $K$-principal bundle $G \to G/K$ by interpreting $Q$ as the homogeneous space $\Gamma\backslash G$. 
\end{prop}
\section{Holomorphic maps from Riemann surfaces}
In this section we suppose that $\pi: P \to M$ is such that $P$ has an integrable almost complex structure $J_\alpha$ associated to a connection $A$ on $P$ and a $\k$--valued 1--form $\alpha$ as in our main construction, and furthermore that $M$ is geodesically complete and connected. By the results of the last section $P$ is then naturally endowed with a locally free, transitive holomorphic right action of the complexified Lie group $G = K^\C$, and with a $\g$--valued holomorphic 1--form $\omega$.

Here we shall be interested in the study of holomorphic curves $f:\Sigma \to P$ in this setting. We shall first see that $f$ is in fact determined by the $\g$--valued 1--form $f^* \omega$ on $\Sigma$. Conversely, we shall ask when a $\g$--valued 1--form $\eta$ on $\Sigma$ gives rise to a holomorphic curve $f: \Sigma \to P$ with $f^*\omega = \eta$. We give a quite complete answer if $\eta$ is of what we shall call `scalar type', but we will show that all such holomorphic maps factor through an elliptic curve, so do not give rise to `interesting' curves of higher genus. This may be compared to Winkelmann's approach  in \cite{Wi4}. 

\subsection{Factorisation properties through elliptic curves}

In the following Lemma we show that an arbitrary smooth pointed map $f: (Y,y_0) \to (P,p_0)$ is determined by the pull--back $f^*\omega$, where $Y$ is a connected smooth manifold and $y_0$ is a point of it. In other words, if $f_1, f_2$ are two such maps with $f_1(y_0) = f_2(y_0)$, and such that $f_1^*\omega = f_2^* \omega$, then in fact $f_1=f_2$. 
\begin{lemma}\label{determination of f}
Let $Y$ be a connected smooth manifold, and let $f: Y \to P$ be smooth and such that $f(y_0)=p_0$, where $y_0$ and $p_0$ are points of $Y$ and $P$ respectively. Then $f$ is determined by $f^*\omega$. 
\end{lemma}
\proof
Let $y_1 \in Y$ be an arbitrary point of $Y$. Let $(y(t))_{t \in [0,1]}$ be a smooth path with $y(0)=y_0$ and $y(1)=y_1$. We claim that there is a unique smooth path $(g(t))_{t \in [0,1]}$ with $g(0) = \id_G$, such that 
\[
  f(y(t)) = p_0 \, g(t) =:p(t)
\]
for all $t \in [0,1]$. In fact, for all $p \in P$ the map $\rho^\C_p: G \to P$, $g \mapsto p g$ is a local diffeomorphism, and so the path $g(t)$ is locally determined. The ambiguity thus remains in the choice of a starting point $g(0) \in \Gamma$, where $\Gamma$ denotes the stabiliser of the point $p_0$ under the $G$-operation. This is fixed by requiring that $g(0)$ is the identity of $G$. Let now $t_0 \in [0,1]$. Differentiating the path $t \mapsto p_0 \, g(t)$ at the point $t_0$ yields the equation
\begin{eqnarray*}
  \dot{p}(t_0) &  = & \dttnull \ p_0 \ g(t) \\
  & = & \dttnull p_0g(t_0) \ g(t_0)^{-1} g(t) \\
  & = & \left(\rho^\C_{p(t_0)} \right) _{*,e} \dtnull g(t_0)^{-1} g(t) \ .
\end{eqnarray*}
By construction of $\omega$ one therefore sees that $(g(t))$ solves the ordinary 
differential equation
\begin{eqnarray}\label{eqgt}
  \dot{g}(t) = g(t) \omega(\dot{p}(t)) \ .
\end{eqnarray}
By the theory of ordinary differential equations there is a unique solution to this equation for a given initial value, and that a solution extends to the domain of definition of $(p(t))$ follows from the lemma in \cite[p.69]{KN1}. But now we simply notice that 
\[
 \omega(\dot{p}(t)) = \omega( \dt f(y(t))) = \omega(f_* \dot{y}(t)) = (f^*\omega)(\dot{y}(t)) \ , 
\]
from which the claim follows. \qed 

The next Lemma is of a similar spirit as the preceeding one, demonstrating the strength of of the holomorphic form $\omega \in \Omega^{1,0}(P,\g)$. Let now $Y$ be a complex manifold. 

\begin{lemma}\label{detection holomorphicity}
Let $f:Y \to P$ be an arbitrary map. Then the following assertions are equivalent:
\begin{enumerate}
	\item The pull-back form $f^*\omega$ is a 1-form of type $(1,0)$, i.e. $f^*\omega \in \Omega^{1,0}(Y;\g)$, and this is a holomorphic form: $\overline{\partial} f^*\omega =0$.  
	\item The pull-back form $f^*\omega$ is a 1-form of type $(1,0)$. 
	\item The map $f$ is holomorphic. \\
\end{enumerate}
\end{lemma}
\proof Clearly $(1)$ implies $(2)$ and $(3)$ implies $(1)$ because $\omega$ is a holomorphic form of type $(1,0)$. Thus it suffices to show that $(2)$ implies $(3)$. We do so locally: We suppose having a complex coordinate chart $\{w^i\}$ on $Y$ and a complex coordinate chart $\{z^j\}$ on $P$ such that $f$ maps the domain of definition of the former to that of the latter. We can write $\omega = \sum \omega_i(z) dz^i$, where $\omega_i(z)$ is a holomorphic function of the coordinates $z=(z_1, z_2, \dots )$. Now read through the charts we have
\begin{align*}
	f^*\omega = \sum_i \omega_i(f(w)) \,  \sum_j \left( \frac{\partial f^i}{\partial w^j} dw^j + \frac{\partial f^i}{\partial \overline{w}^j} d\overline{w}^j\right)  \ . 
\end{align*}
The assumption that $f^*\omega$ is of type $(1,0)$ now implies that for all $j$ we have 
\begin{align*}
	0 = \omega\left(\frac{\partial f}{\partial \overline{w}^j} \right) = \sum_i \omega_i(f(w)) \,   \frac{\partial f^i}{\partial \overline{w}^j}  
\end{align*}
at all points of the coordinate chart $\{ w^{i} \}$.
However, pointwise $\omega: TP \to \g$ is an isomorphism, and  so the last equation implies that 
\[
\frac{\partial f}{\partial \overline{w}^j}  = 0 
\]
for all $j$, thus implying that $f$ is holomorphic. 
\qed

As the Proof of Lemma \ref{determination of f} relied on the uniqueness of the solution of an ordinary differential equation, one might ask whether the following converse of Lemma \ref{determination of f} and Lemma \ref{detection holomorphicity} is true: Given a (holomorphic) form $\eta \in \Omega^1(Y,\g)$ on $Y$, is there a (holomorphic) map $f: Y \to P$ such that $f^* \omega = \eta$? Inspired by the proof of Lemma \ref{determination of f} the idea would certainly lie in trying to define $f$ in the following way: Again, we assume having a chosen point $z_0 \in Y$ and a point $p_0 \in P$, and we will require that $f(z_0) = p_0$ (this is of course no real restriction as the $G$--action is transitive on $P$). Now let $z \in Y$ be an arbitrary point, and let $\tau_z: [0,1] \to Y$ be a smooth path from $z_0$ to $z$. We would seek a solution $g_{\tau_z} : [0,1] \to G$ of the differential equation
\begin{equation} \label{diff eq g}
	\dot{g}_{\tau_z}(t) = g_{\tau_z}(t) \, \eta(\dot{\tau}_{z}(t)) \ ,
\end{equation}
with the initial condition $g_{\tau_z}(0)=1$, 
and we would define 
\begin{equation} \label{defining f}
f(z):= p_0 \, g_{\tau_z}(1) \ .
\end{equation}
This must be independent of the chosen path $\tau_z$ from $z_0$ to $z$. This will be the case if for any {\em closed} loop $\gamma: [0,1] \to Y$ with start- and endpoint $y_0$ the solution $g_{\gamma}$ of $\dot{g}_{\gamma} = g_{\gamma} \, \eta(\dot{\gamma})$ is such that the endpoint $g_{\gamma}(1)$ lies in the stabiliser $\Gamma$ of the $G$--action on $P$ at $p_0$. Indeed, if this is the case, then for any two paths $\tau_z$ and  $\tilde{\tau}_z$ from $z_0$ to $z$ the concatenation $\tilde{\tau}^{-1} * \tau$ will be a closed loop based at $z_0$, and for the endpoint $g_{\tilde{\tau}^{-1} * \tau}(2)$ of the solution to the above differential equation we will have
\[
   g_{\tilde{\tau}^{-1} * \tau}(2) = g_{\tau_z}(1) g_{\tilde{\tau}_z} (1)^{-1}  \in \Gamma \ ,
\]
and so the definition $f(z):= p_0 \, g_{\tau_z}(1)$ is independent of the chosen path \footnote{We have been a little abusive here with parametrisation of the paths - of course, the endpoints $g_{\tau_z}(1)$ are independent of the parametrisation chosen for the path $\tau_z$, and so it doesn't really matter that the concatenation is a priori only piecewise smooth.}. 

We give an affirmative answer to this question in the case where the holomorphic form $\eta \in \Omega^1(Y;\g)$ is closed and of a special type: 
\begin{prop}
  Let $\eta$ be a closed holomorphic form with value in $\g$ which is of ``scalar type'' in the following sense: There is a holomorphic differential form $\zeta$ and an element $Z \in \g$ such that 
$$
\eta = Z \cdot \zeta \ .
$$
1. There exists a holomorphic map $f: Y \nach P$ with $f^* \omega = \eta$ and $f(z_0) = p_0$ 
if and only if 
\begin{eqnarray}\label{expzetah1}
\exp(\displaystyle Z \cdot  [\zeta]_{dR} (H_1(Y;\Z))) \subseteq \Gamma \ .
\end{eqnarray}
2. If this is the case, then each such holomorphic map factorises through an elliptic curve. This means that we have a commutative diagram
\begin{eqnarray}
\begin{diagram}
\node{Y} \arrow[1]{e,t}{f} \arrow[1]{s,l}{\bar{f}} \node[1]{P}  \\
\node{\Sigma_0} \arrow[1]{ne,r}{} 
\end{diagram} 
\end{eqnarray}
where all maps are holomorphic, and where $\Sigma_0$ is an elliptic 
curve.
\end{prop}
\begin{remark}
1. We would like to point out that if $Y$ is a Riemann surface $\Sigma$, then we do not need to ask that the holomorphic form $\eta$ is closed. It will be so automatically. 
\\
2. In the above formula (\ref{expzetah1}) the form $[\zeta]_{dR}$ is interpreted via the isomorphism \\ $H^1_{dR}(Y,\C) \cong \Hom_\Z (H_1(Y;\Z),\C)$.
\end{remark}
\proof
Let $z \in \Sigma$ be an arbitrary point, and let $[0,1] \to Y$, $t \mapsto \tau_z(t)$ be a path from $z_0$ to $z$. 
The fact that $\eta$ is of scalar form implies that the equation
\begin{eqnarray}\label{eq_g_nonclosed}
 \dot{g}_\tau (t) = g_\tau (t) \ \eta(\dot{y}(t))
\end{eqnarray}
is easily integrated by the defining relation of the exponential map \cite{War}. For  a curve $\tau_z(t), t \in [0,1]$ the Lie group element $g_{\tau_z}(t)$ is simply given by
\begin{equation}\label{g with exponential}
g_{\tau_z}(t) = \exp( Z \cdot  \int_0^t \zeta(\dot{\tau}_z (s)) ds ) \ .
\end{equation}
We see in particular, that for a closed curve $\tau$ the integral is just the cohomology-homology pairing:
$$
\int_0^1 \zeta(\dot{\tau} (s)) ds = [\zeta]_{dR} ([\tau]) \ .
$$

To prove the first statement, let first $f: Y \to P$ be a holomorphic (pointed) map with $f^* \omega = \eta$. As we have seen in the proof of Lemma \ref{determination of f}, $f(z)$ has to satisfy the equation (\ref{defining f}) with $g_{\tau_z}(1)$ satisfying the differential equation (\ref{diff eq g}).
In particular, we need  $g_\gamma (1) \in \Gamma$ for any closed curve $\gamma$ in order to have $f$ well defined. This must also be true for $\gamma$ an arbitrary cycle. By the above formula (\ref{g with exponential}) and the interpretation of the integral in terms of the pairing between homology and cohomology, we get (\ref{expzetah1}).
\\

For the contrary implication, assume that the formula (\ref{expzetah1}) holds. As explained above, this implies that the definition  
$$
f(z) := p_0 \ g_{\tau_z}(1) \ ,
$$
is well defined, where $\tau_z : [0,1] \nach Y$ is an arbitrary curve from $z_0$ to $z$.
By construction we have $f^* \omega = \eta$. Now Lemma \ref{detection holomorphicity} implies that $f$ is holomorphic. 
\\

To prove the second assertion, we have seen that the assumption implies that we have the formula
\begin{eqnarray*}
  f(z) = p_0 \ \exp ( Z \cdot \int_0^1 \zeta(\dot{\tau}_z (t)) dt)
\end{eqnarray*}
for the map $f$. Of course, this factors
%where $\tau_z$ is an arbitrary path from $z_0$ to $z$. As we have seen the 
%choice of another curve $\tau_z$ from $z_0$ to $z$ will make differ the exponent by an element of the abelian group $Z \ [\zeta]_{dR}(H_1(Y;\Z))$. But these elements are mapped via the exponential to the stabilisator $\Gamma$, so that $f$ factors 
in the following way,
\begin{eqnarray}
\begin{diagram}
\node{Y} \arrow[1]{e,t}{f} \arrow[1]{s,l}{\tilde{f}} \node[1]{P \ ,}  \\
\node{\Cquotzetah} \arrow[1]{ne,r}{} 
\end{diagram} 
\end{eqnarray}
Note that it doesn't in general even make sense to ask whether $\tilde{f}$ is a holomorphic (or differentiable) map, since in general the quotient $\Cquotzetah$ does not even need to have the structure of a manifold. The group $\zeth$ may in general be quite complicated, for instance, it can even be dense in $\C$. However, the fact that $\Gamma$ is a closed subgroup of $G$ implies even that $f$ factorises through the quotient of $\C$ by the closure of the group $\zeth$. The group $\overline{\zeth}$ being an abstract closed subgroup of the Lie group $\C$, it must be a Lie subgroup by the Theorem of Chevalley. Thus $f$ factorises as
\begin{eqnarray}
\begin{diagram}
\node{Y} \arrow[1]{e,t}{f} \arrow[1]{s,l}{\bar{f}} \node[1]{P \ ,}  \\
\node{\Cquotclzetah} \arrow[1]{ne,r}{p_0 \exp(Z \cdot \, \_)} 
\end{diagram} 
\end{eqnarray}
where $\bar{f}$ is at least differentiable, since the projection $\C \to\Cquotclzetah$ endowes this quotient with the structure of a smooth manifold. For the time being, we shall assume that $\overline{\zeth}$ is a discrete subgroup of $\C$, so that $\Cquotclzetah$ is a manifold of real dimension $2$, naturally endowed with a complex structure. In this case $\bar{f}$ is holomorphic. To see this, we can consider a local lift of $\bar{f}$ to $\C$: Let $(U,z)$ be a coordinate chart\footnote{to keep notations simple, we make no distinction in notation between preimage and image of a chart} around a point $z_1 \in Y$. We may assume that $U$ is the open unit ball in $\C^{\dim Y}$, and that $z_1 = 0$. As path $\tau_z$ we will then choose a path $\tau_{z_1}$ followed by the path $[0,1] \to U$, $t \mapsto t \, z$, which is the traight line from $z_1$ to $z$ in the coordinate chart. Writing $\zeta$ as $\zeta = \sum \zeta_{i} \, dz^i$, with holomorphic $\g$--valued functions $\zeta_i$, we obtain the following formula for a local lift of $\bar{f}$:
\begin{align*}
z \mapsto \int_0^1 \zeta(\dot{\tau}_z (t)) \, dt 
		= \int_0^1 \zeta_i(t z) \cdot z^i \,  dt \, .
\end{align*}
As integration commutes with the differential operators $\partial/ \partial \bar{z}^{j}$ (by using, for instance, Lebesgue's theorem of dominated convergence), we see that this lift is indeed holomorphic, and so $\bar{f}$ is holomorphic, too.  

We will finally discuss the possible cases for $\clzeth$, starting with positive dimension. If we have $\dim_\R \clzeth = 2$, then $f$ factorises over a point, thus is constant (and therefore factorises obviously over an elliptic curve). If $\dim_\R \clzeth = 1$, then $f$ must be constant, too. If it were not, $f$ would be of rank at most 1. But a non--constant holomorphic map has a non-vanishing derivative somewhere, and this has real dimension at least 2.

We make three distinction in the final case $\dim_\R \clzeth = 0$. 
\begin{enumerate}[label=(\alph*)]
\item If $\clzeth = \{0\}$, then $\bar{f}$ becomes a holomorphic map to $\C$, hence has to be constant since $Y$ is compact (and thus factors through an elliptic curve).
\item If $\clzeth = \left<b_1\right>_\Z$ with $b_1 \in \C^*$ then $\Cquotclzetah \cong \C^* \subseteq \C$. Thus $\bar{f}$ is a holomorphic map to $\C$, hence has to be constant as $Y$ is compact, and therefore factors through an elliptic curve.
\item If $\clzeth = \left<b_1, b_2\right>_\Z$ where $b_1, b_2$ are $\R$-linearly independent in $\C$, hen $f$ factors through $\C / \left<b_1,b_2\right>$ which is indeed an elliptic curve.
\end{enumerate}
\qed

\begin{remark}
\begin{enumerate}
\item[1.] An alternative viewpiont of this result is the following: If one has given a holomorphic curve $f: \Sigma \to P$ from a Riemann surface $\Sigma$ of genus $g \geq 2$, and one wants to be sure that this does not factor through an elliptic curve, then it is necessary that the holomorphic $1$--form $f^* \omega$ is not of scalar type. 
\item[2.] As a corollary from the proof we get that there are no (non-constant) rational curves $f: \mathbb{C}P^1 \to P$ with $f^* \omega$ of scalar type. More generally, if the structure group $K$ is abelian, then there are no (non-constant) rational curves $f: \C P^1 \to P$ at all. Of course, given that $\C P^1$ has no non-trivial holomorphic forms, this also follows from the Lemmata \ref{determination of f} and \ref{detection holomorphicity} above.
\end{enumerate}
\end{remark}
\qed
\subsection{Another factorisation property}
Again, let $K$ be a compact Lie group with complexification $G = K^\C$, and let $\g$ be the complex Lie algebra of $G$. 
Let $\Sigma$ be a Riemann surface, and let 
$\eta \in \Omega^1(\Sigma,\g)$ be a non-trivial holomorphic 1-form with value in $\g$. To $\eta$ we can associate a map $\phi_\eta$ to the projectivisation $P(\g)$ (a complex projective space of dimension $\dim_\C(\g)$)  as follows. 

Let $(U,z)$ be a complex coordinate chart of $\Sigma$. Then on $U$ the form $\eta$ can be written as
$$
\eta = \mu \, dz \ , 
$$
where $\mu: U \to \g$ is a holomorphic $\g$-valued function on $U$. Now if $V \subseteq U$ is a domain where $\mu$ has no zeros, then $\phi_\eta$ is simply defined by
\[
	\phi_\eta(z) := [\mu(z)] \in P(\g) \ .
\]
Of course the element $\mu(z) \in \g$ associated to $\eta$ in the above manner depends on the coordinate chart, but its associated element in the projectivisation $P(\g)$ doesn't. 
Suppose now that the holomorphic map $z \mapsto \mu(z)$ vanishes at $z_0$. Since $\eta$ is non-trivial there must be a positive integer $k \in \N$ such that the expression
\[
\tilde{\mu}(z) := \frac{\mu(z)}{(z-z_0)^k}
\]
defines a holomorphic map in a neighbourhood of $z_0$ which is non-vanishing in $z_0$. Now define $\phi_\eta(z) := [\tilde{\mu}(z)]$ in a neighbourhood of $z_0$. Away from $z_0$ the two definitions coincide, so the latter extends the former across $z_0$. As this definition behaves well under change of coordinates, patching this definition yields a well defined map $\phi_\eta : \Sigma \to P(\g)$.

Let us now suppose that we have a holomorphic curve $f: \Sigma\nach P$ with induced map $\bar{f}: \Sigma\nach M$ defined by $\bar{f}= \pi \circ f$. We call $\bar{f}$ {\em conformal} with respect to the canonical almost complex structure $J_0$ of $\Sigma$ if\footnote{$\bar{f}^*g$ is not required to be a metric on $Y$. It is so if and only if $\bar{f}$ is an immersion} 
$$
  \bar{f}^* g \ (u,  J_0 u) = 0 
$$
holds for all $u \in T\Sigma$. We shall point out that if $\bar{f}$ is an immersion then the condition on $\bar{f}$ to be conformal means that the complex structure induced by $\bar{f}^*g$ (and the given orientation) coincides with the initial complex structure $J_0$ on $\Sigma$. Here $g$ is supposed to be the metric $g_\alpha$ introduced above -- the one that makes  the map $\val$ from above an isometry between the associated vector bundle $P \times_{\ad} \k$ and the tangent bundle $TM$. 

We get the following interesting algebraic geometrical fact relating the conformality of $\bar{f}$ to the factorisation of $f$ through a quadric. 

\begin{prop}
Suppose the non-constant holomorphic map $f:\Sigma\nach P$ induces a conformal map $ \bar{f}: \Sigma\nach M$. Then the canonically induced holomorphic map $\phi_\eta: \Sigma\nach P(\g)$, associated to the holomorphic form $\eta:=f^*\omega$, factorises through a smooth quadric $Q$ in $P(\g)$. In other words, there is a complex codimension-$1$-submanifold $Q$ in $P(\g)$, defined by a polynomial of degree 2 on  $\g$. 
\begin{eqnarray*}
\begin{diagram}
\node{\Sigma} \arrow[2]{e,t}{\phi_\eta} \arrow[1]{se,l}{\psi_\eta} \node[2]{P(\g) } \\
\node{} \node[1]{Q} \arrow[1]{ne,r}{}
\end{diagram} 
\end{eqnarray*}
\end{prop}
\begin{remark}
This is particularly intersting in the case where $K$ is $SU(2)$ or $SO(3)$. In this case any non-constant holomorphic map $f: \Sigma \to P$ which induces a  conformal map $\bar{f}: \Sigma \to M$ to the base is such that the canonical map $\phi_\eta: \Sigma \to P^2(\C) \cong P(\so(3;C))$ associated to $\eta= f^* \omega$ factors through a quadric $Q$ in $P^2(\C)$,  or in other words, through a rational curve. 
\end{remark}

\proof
%If we have a holomorphic map $f:\Sigma\nach P$, then we will consider the case where $f$, composed with the projection $\pi$, becomes an immersion $ \bar{f}: \Sigma\nach M$. In that case the holomorphic form $\eta:=f^*\omega$ with value in the Lie algebra $\so(3,\C)$ induces a canonical holomorphic map $\phi_\eta: \Sigma\nach P^2\C$. In the particular case where the complex structure $J_0$ on $Y$ is equal to the one which is induced via $\bar{f}$ from the metric $g_\alpha$ on $M$, the following happens: $\phi_\eta$ factorises through a quadric $Q$ in $P^2 \C$ which is biholomorphic to $P^1 \C$.
There is a unique (sesquilinear) scalar product $\langle \, . \, , \, . \, \rangle_{\g}$ on $\g$ which coincides with the inner product on $\k$ on real elements. More precisely, this inner product is given as the tensor product of the inner product on $\k$ with the sesquilinear form $(z,w) \mapsto \bar{z} w$ on $\C$, where $\g = \k \tensor_{\R} \C$. 

We will have to relate the fact that $\bar{f}$ is conformal with the map $\phi_\eta$ associated to $\eta= f^* \omega$. To do so, we will use the fact that $f$ is holomorphic in order to reexpress the pull-back $\bar{f}^*g$. 
That $f$ is holomorphic means that its derivative $f_*$ commutes with the corresponding almost complex structures: $\Ja \circ f_* = f_* \circ J_0$. By the very definition of $\Ja$ we get
\begin{eqnarray*}
\Ja \circ f_* (u) & = & \alpha^\# (f_* u) 
- \alpha^{-1}_{f(z)} (\omega_A(f_* u))
\end{eqnarray*}
for an arbitrary element $u \in T_z\Sigma$ on the one hand, as well as 
\begin{eqnarray*}
  f_* \circ J_0 (u) = (\omega_A(f_*(J_0 u)))^\#_{f(z)} 
  + \alpha^{-1}_{f(z)} (\alpha(f_* (J_0 u)))
\end{eqnarray*}
on the other hand. Comparison of the vertical and the horizontal parts gives now the following fact: $f$ is holomorphic if and only if
\begin{eqnarray*}
\alpha(f_* (u)) = (f^* \omega_A)(J_0 u)  
\end{eqnarray*}
for all $u \in T\Sigma$. 
With the fact that $\val$ is an isometry we therefore get
\begin{eqnarray}\label{metrique_hol}
\bar{f}^*g_\alpha \ (u , v) = 
\langle (f^* \omega_A)(J_0 u),  (f^* \omega_A)(J_0 v) \rangle_{\k} \ .
\end{eqnarray}
\\

The $\k$-valued 1-form $\omega_A$ is just the real part of the holomorphic $\g$-valued 1-form $\omega$. Thus $f^* \omega_A = f^* \Re \  \omega = \Re \ f^* \omega$. Note that we have $\omega \circ \Ja = i \omega$. For $\eta = f^* \omega$ the above equation (\ref{metrique_hol}) yields the following expression for the pull-back of $g$ via $\bar{f}$: 
\begin{equation}\label{pull-back of metric}
\begin{split}
\bar{f}^*g \ ( u , v) 
& =  \langle \, \Re \  \eta (J_0 u), \Re \  \eta (J_0 v) \, \rangle_{\k} \\
& =  \langle \, \Im \  \eta (u), \Im \  \eta (v) \, \rangle_{\k} \\
& =  - \frac{1}{4} \langle \,  \eta(u)- \overline{\eta(u)} , 
                       \eta(v) - \overline{\eta(v)} \,  \rangle_{\g} \\
& =  - \frac{1}{2} \, ( \Re \  \langle\, \eta(u),\eta(v)\rangle_{\g} - \Re \  \langle\eta(u), \overline{\eta(v)}\rangle_{\g}) \ .
\end{split}
\end{equation}
As $\bar{f}$ is supposed to be conformal by hypothesis we get the equation 
\begin{equation}\label{quadric 1}
\begin{split}
 0 =   \bar{f}^*g \  (u,J_0 u) & =  - \frac{1}{2} (\Re \  \langle \eta(u),i \eta(u)\rangle_{\g} - \Re \  \langle \eta(u), \overline{i \eta(u)} \rangle_{\g}) \\
  & =  \frac{1}{2} (\Im \  \langle \eta(u), \eta(u) \rangle_{\g} - \Im \  \langle \eta(u), \overline{\eta(u)} \rangle_{\g}) \\
  & =  - \frac{1}{2}\, \Im \, \langle\eta(u), \overline{\eta(u)} \rangle_{\g} 
\end{split}
\end{equation}
for all $u \in T\Sigma$. 

In order to see that $\phi_\eta$ factors through a smooth quadric, let us choose an orthonormal basis $(X_1, \dots, X_N)$ of $\k$, yielding also an orthonormal basis for the sesquilinear extension on $\g$. This also yields an identification $P(\g) \cong P^{N-1} \C$. 

Let $(U,z)$ be a chart of holomorphic coordinates of $\Sigma$, and let $\eta = \mu \, dz$, where $\mu$ is a $\g$-valued holomorphic function. We then have $\phi_\eta(z) = [\mu_1(z), \dots, \mu_N(z)]$ after the identification induced by the basis $(X_i)$ above, where we have expressed $\mu(z) = \sum_{i=1}^{N} \mu_i(z) \, X_i$. With complex coordinates $z = x + iy$ we have
\[
	\eta\left(\frac{\partial}{\partial x}\right) = \mu \ ,
\]
Therefore, the condition of $\bar{f}$ to be conformal implies, via the above equation (\ref{quadric 1}), that on the chart $(U,z)$ we have in particular 
\begin{equation*}
	0 = \Im \, \langle \, \eta\left(\frac{\partial}{\partial x}\right) , \overline{\eta\left(\frac{\partial}{\partial x}\right)} \, \rangle_{\g}
		= \Im \sum_{i=1}^{N} (\mu_i(z))^2 \  
\end{equation*}
at any point in $U$. But as the function $z \mapsto \sum (\mu_i(z))^2$ is holomorphic the vanishing of its imaginary part implies the vanishing of its real part either. Therefore we have 
\begin{equation*}
   \sum_{i=1}^{N} (\mu_i(z))^2 = 0 \ 
\end{equation*}
at any point of $U$. A similar discussion applies at points $z_0$ where $\mu$ vanishes. Also, this holds for any chart of a compatible holomorphic atlas of $\Sigma$, and so $\phi_\eta$ indeed factors through a smooth quadric.

%Thus the first compatibility condition $\bar{f}^*g \ (u,J_0 u) = 0$ becomes 
%\begin{eqnarray}\label{first_compatibility}
%  \Im \  \langle \eta(u), \overline{\eta(u)} \rangle \ =
%\Im \  \left( \sum \ (\eta^i(u))^2 \right) = 0 \ .
%\end{eqnarray}
%For the second compatibility condition we calculate 
%\begin{eqnarray*}
%\bar{f}^*g \  (u,u) & = & \frac{1}{2} (\Re \  \langle\eta(u),\eta(u)\rangle_{\so(3,\C)} - \Re \  \langle\eta(u), \overline{ \eta(u)}\rangle_{\so(3,\C)})
%\end{eqnarray*}
%and also
%\begin{eqnarray*}
%\bar{f}^*g \  (J_0 u,J_0 u) & = & \frac{1}{2} (\Re \  \langle\eta(u),\eta(u)\rangle_{\so(3,\C)} - \Re \  \langle \eta(u), -i \overline{ \eta(u)}\rangle_{\so(3,\C)}) \\
%& =  & \frac{1}{2} (\Re \  \langle\eta(u),\eta(u)\rangle_{\so(3,\C)} + \Re \  \langle \eta(u), \overline{ \eta(u)}\rangle_{\so(3,\C)}) \ .
%\end{eqnarray*}
%Thus the equailty implies 
%\begin{eqnarray}
%  \Re \  ( \langle \eta(u), \overline{ \eta(u)}\rangle_{\so(3,\C)} ) = \Re \  \left( \sum (\eta^i(u))^2) \right) = 0
%\end{eqnarray}
%which gives the one and only compatibility equation 
%\begin{eqnarray*}
%  \sum_{i=1}^{3} (\eta^i(u))^2 = 0 \quad \quad \forall u \in T\Sigma.
%\end{eqnarray*}

\qed

\end{document}